\documentclass{amsart}

\usepackage{verbatim}
\usepackage{xcolor}
\usepackage{tikz}
\usepackage{tikz-cd}
\usepackage{arydshln}
\usepackage{caption}
\usepackage{subcaption}
\usepackage{graphicx}
\tikzset{
	symbol/.style={
		draw=none,
		every to/.append style={
			edge node={node [sloped, allow upside down, auto=false]{$#1$}}}
	}
}

\usepackage{amsmath}
\usepackage{bbm}
\usepackage{mathtools}

\usepackage{hyperref}
\usepackage{times,amsfonts,amsmath,amstext,amsbsy,amssymb,
amsopn,amsthm,upref,eucal,mathrsfs}
\usepackage[T1]{fontenc}

\newtheorem{theorem}{Theorem}[section]
\newtheorem{proposition}[theorem]{Proposition}
\newtheorem{lemma}[theorem]{Lemma}

\theoremstyle{definition}
\newtheorem{definition}[theorem]{Definition}
\newtheorem{example}[theorem]{Example}

\theoremstyle{remark}

\numberwithin{equation}{section}

\newcommand{\opnorm}[1]{\|#1\|_{\mathrm{op}}}

\newcommand{\field}[1]{\mathbb{#1}}
\newcommand{\R}{\field{R}}
\newcommand{\N}{\field{N}}

\newcommand{\Z}{\field{Z}}

\newcommand{\GL}{\mathrm{GL}}
\newcommand{\Pb}{\mathbf{P}}
\newcommand{\E}{\mathbf{E}}
\newcommand{\Cal}{\mathcal}
\newcommand{\A}{\Cal{A}}
\newcommand{\V}{\Cal{V}}
\newcommand{\B}{\mathscr{B}}

\begin{document}

\title[From limit theorems to mixing limit theorems]{From limit theorems  to mixing limit theorems}


%    author one information
\author{Jianshu Jiang}
%\address{Department of Mathematics, Stanford University, 450 Serra Mall
%	Building 380, Stanford, CA 94305-2125, USA}
%\curraddr{}
%\email{farana@stanford.edu}
%\thanks{}

\email{jianshu.jiang@cyu.fr}

\address{Laboratoire AGM Batîment E, 5ème étage 2 av. Adolphe Chauvin
95302 Cergy-Pontoise Cedex, France}

%    author two information
%\author{}
%\address{}
%\curraddr{}
%\email{}
%\thanks{}

%\subjclass[2000]{Primary }
%    The 2010 edition of the Mathematics Subject Classification is
%    now available.  If you are citing a classification from the
%    new scheme, use the following input coding instead.
%\subjclass[2010]{Primary }

\date{}

\begin{abstract}
    Following the work of Gouëzel, we give an alternative proof, under weaker assumptions, of a result of Arana-Herrera and Forni, which provides a general method for upgrading limit theorems for cocycles over dynamical systems to mixing limit theorems.
\end{abstract}

\maketitle


\thispagestyle{empty}

\section{Introduction}

In ergodic theory, several fundamental notions are employed to capture the intrinsic randomness of a deterministic dynamical system. At the most basic level, ergodicity ensures that for almost every point, its time averages along orbit converge to the spatial average. To obtain a finer statistical description, distributional limit theorems—most notably the central limit theorem (CLT)—are used to characterize the asymptotic fluctuations of time averages around their spatial mean. From a different perspective, mixing, a property strictly stronger than ergodicity, guarantees the asymptotic independence between the initial state of a trajectory and its position at large forward times.

While classical limit theorems focus on Birkhoff sums of scalar observables, great effort has been devoted to understanding non-abelian analogues, namely linear cocycles over dynamical systems. In this setting, the natural observables are the logarithmic growth of a vector or of the operator norm along a typical orbit. Central limit theorems of this kind have been established in a variety of settings: Benoist, Quint\cite{benoist_quint2016} and more recently Grama, Quint, Xiao~\cite{grama2026locallimittheoremoperator} treated products of independent random matrices, Furman and Kozma \cite{furman_kozma2021} studied matrix cocycles driven by Markov measures on subshifts of finite type, and Al-Saqban and Forni \cite{alsaqban_forni2025} proved central limit theorems for the Kontsevich--Zorich cocycle over the Teichm\"uller flow.

Whenever a dynamical system and an associated observable or cocycle satisfy a distributional limit theorem (DLT), one is naturally led to investigate the statistical ``robustness'' of the limiting law under different observational conditions. In many physical or geometric situations, the phase space comes equipped with an \emph{a priori} reference probability measure $\nu$ (for instance, the Lebesgue or volume measure) which may differ from the dynamically invariant measure $\mu$. It is then natural to ask whether distributional convergence holds with respect to the physical reference measure $\nu$ rather than only with respect to $\mu$. This directly motivates our first question:

\medskip
\noindent\textbf{Question 1:} \emph{If a DLT is established under the invariant measure $\mu$, can we deduce the same asymptotic distributional behavior under another reference measure $\nu$?}
\medskip

In the scalar setting of Birkhoff sums (or additive real-valued cocycles), this question is answered affirmatively by a classical theorem of Eagleson \cite{eagleson1976some}, which asserts that for an ergodic probability measure-preserving system, any distributional limit law remains unchanged under an absolutely continuous change of the underlying probability measure ($\nu \ll \mu$); see also Zweim\"uller \cite{Zweimuller2007} for a different proof. In particular, choosing $\nu$ to be the normalized restriction of $\mu$ to a measurable set of positive measure shows that conditioning on the initial position of an orbit does not alter the limiting distribution. Beyond its intrinsic interest as a one-sided conditioning principle, Eagleson's theorem plays a crucial technical role in deducing limit theorems for a dynamical system from those of an induced transformation: a key step in such arguments is to replace the invariant measure on the inducing set by a measure weighted by the return time while preserving the limit law; see, for instance, \cite{Melbourne2004,Chazottes2007}.

Going one step further, one may ask what happens if we condition not only on the initial state, but simultaneously on the initial and terminal states—or even on multiple states—along the orbit. This leads to our second question:

\medskip
\noindent\textbf{Question 2:} \emph{Does the limiting distribution remain unchanged if we prescribe both the initial and final positions of the orbit segments?}
\medskip

Mathematically, simultaneous conditioning at the initial and terminal times is considerably more delicate than the one-sided problem because the terminal condition moves with the orbit length. Building on Zweim\"uller's approach \cite{Zweimuller2007}, Gou\"ezel \cite{GOUËZEL_2020} resolved this two-sided problem for Birkhoff sums, showing that if the base transformation is mixing, distributional limit theorems remain valid after prescribing both endpoints. 

More recently, Arana-Herrera and Forni \cite{arana_herrera_forni2024} developed a higher-dimensional analogue for linear cocycles and applied it to the Kontsevich--Zorich cocycle over the Teichm\"uller flow. Assuming that the base dynamics is mixing, together with suitable hyperbolicity assumptions on the base and regularity conditions on the cocycle, they showed that distributional limit theorems for logarithmic cocycle growth can be upgraded to mixing limit theorems. In forthcoming work of Arana-Herrera and Honaryar \cite{AranaHerrera2023}, these results are further applied to establish laws of large numbers and central limit theorems for the arithmetic and homological complexity of long simple closed geodesics on negatively curved surfaces.

The main contribution of the present note is to re-prove and extend the Eagleson-type and mixing limit theorems of Arana-Herrera and Forni \cite{arana_herrera_forni2024} under substantially weaker hypotheses. Specifically, we provide a direct criterion for finite-dimensional linear cocycles when projectivized directions are averaged against the standard rotation-invariant probability measure on projective space. Our formulation requires \emph{only} the mixing property of the base transformation and the log-integrability of the cocycle—dispensing entirely with any hyperbolicity assumption on the base or any regularity condition on the cocycle beyond integrability.

Thus, for any mixing example among the cocycles cited above \cite{benoist_quint2016,furman_kozma2021,alsaqban_forni2025}, our main theorem provides an immediate route from an unconditioned DLT to a limit law conditioned on both endpoints. Beyond the geometric applications in \cite{arana_herrera_forni2024,AranaHerrera2023}, we expect mixing limit theorems of this type to serve as a tool as the statistical theory of linear cocycles continues to develop.

The remainder of this note is organized as follows. Section~2 introduces the precise setup and states our main results (Theorems~\ref{thm:eagleson-cocycle} and~\ref{thm:main}). Section~3 is devoted to the proofs, proceeding via the reduction to the operator norm, asymptotic invariance under fixed shifts, and a mixing averaging estimate.

\medskip\noindent\textbf{Acknowledgements.} The author would like to thank his supervisor Prof. Giovanni Forni for introducing
him to the subject of mixing limit theorems. The author would also like to thank Dr. Pouya Honaryar for very enlightening conversations on the subject of this paper.

\medskip\noindent\textbf{AI Disclosure.} The main results of the presenting paper was proved without Generative AI tools.

\section{Statement of the results}

Let $(X,\mathscr B,\mu)$ be a Borel probability space and $T:X \to X$ be an invertible probability measure-preserving transformation. Fix $m \in \N$. By an $m$-dimensional linear cocycle over $T$ we mean a measurable map $C : X \times \Z \to \GL(m, \R)$ satisfying the cocycle identity
\[
C(x,r+s)=C(T^rx,s)\cdot C(x,r)\quad \text{for } \mu\text{-almost every } x\in X \text{ and every } r,s\in\Z.
\]

Denote by $\langle \cdot, \cdot \rangle$ the standard inner product on $\R^m$, by $\|\cdot\|$ the corresponding Euclidean norm, by $\opnorm{\cdot}$ the induced operator norm on $\GL(m,\R)$, and by $\eta$ the rotation-invariant probability measure on the projective space $\mathbb{P}\R^m$. The projectivized bundle $X \times \mathbb{P}\R^m$ is endowed with the product probability measure $\mu \otimes \eta$.

Given $(x,v) \in X \times (\R^m \setminus \{0\})$ (or $X \times \mathbb{P}\R^m$) and $N \in \mathbb{N}$, we consider the logarithmic cocycle growth functions:
\begin{align*}
    \sigma(x, v, N) &\coloneqq \log \frac{\|C(x, N)v\|}{\|v\|} \in \R, \\
    \sigma(x, N) &\coloneqq \log \sup_{v \in \R^m \setminus \{0\}} \frac{\|C(x, N)v\|}{\|v\|}=\log \opnorm{C(x,N)} \in \R.
\end{align*}
We say that the cocycle $C$ is \emph{log-integrable} if the maps
\[
    x \mapsto \sigma^+(x,1):=\max\{0, \sigma(x, 1)\}, \quad x \mapsto \sigma^+(x,-1):=\max\{0, \sigma(x, -1)\}
\]
belong to $L^1(X, \mu)$.

\begin{definition}
    We say that the cocycle $C$ satisfies a \emph{(spatial) distributional limit theorem} (DLT) on $(X \times \mathbb{P}\R^m, \mu \otimes \eta)$ if there exist an \emph{averaging sequence} $\A := (A_N)_{N\in\mathbb{N}}$ of real numbers, a \emph{normalizing sequence} $\V := (V_N)_{N\in\mathbb{N}}$ of positive real numbers with $V_N \to \infty$ as $N \to \infty$, and a real-valued random variable $S$ (the \emph{limiting distribution}), such that the random variables
\[
    S_N(x, v) \coloneqq \frac{\sigma(x, v, N) - A_N}{V_N} \quad \text{on } (X \times \mathbb{P}\R^m, \mu \otimes \eta)
\]
converge in distribution to $S$ as $N \to \infty$; that is, for every interval $(a, b) \subseteq \R$ with $\Pb(S \in \{a, b\}) = 0$,
\[
    \lim_{N \to \infty} (\mu \otimes \eta)\big(\{(x, v) \in X \times \mathbb{P}\R^m \mid S_N(x, v) \in (a, b)\}\big) = \Pb(S \in (a, b)).
\]
\end{definition}

\begin{example}
    By the Oseledets multiplicative ergodic theorem, if $T$ is ergodic, every log-integrable cocycle $C : X \times \Z \to \GL(m, \R)$ with top Lyapunov exponent $\lambda \in \R$ satisfies a DLT with averaging sequence $A_N = \lambda N$, normalizing sequence $V_N = N$, and limiting distribution $S \equiv 0$ (the Dirac mass at $0$).
\end{example}

\begin{example}
    Another central example in hyperbolic dynamics is the central limit theorem (CLT), where $A_N = \lambda N$, $V_N = \sqrt{N}$, and the limiting random variable $S$ follows a centered normal distribution $\mathcal{N}(0, \sigma^2)$.
\end{example}

Our first result addresses \textbf{Question~1} in the Introduction, establishing an Eagleson-type theorem for log-integrable cocycles over ergodic systems:

\begin{theorem}[Eagleson's theorem for cocycles]
\label{thm:eagleson-cocycle}
Let $(X,\mathscr B,\mu,T)$ be an ergodic invertible probability measure-preserving system, and let $C:X\times\mathbb Z\to\mathrm{GL}(m,\mathbb R)$ be a log-integrable cocycle. Suppose that $S_N$ converges in distribution to $S$ on $(X\times\mathbb P\mathbb R^m,\mu\otimes\eta)$. Then, for every probability measure $\nu$ on $(X,\mathscr B)$ with $\nu\ll\mu$, the sequence $S_N$ also converges in distribution to $S$ on $(X\times\mathbb P\mathbb R^m,\nu\otimes\eta)$; that is, for every interval $(a,b) \subseteq \R$ such that $\Pb(S\in\{a,b\})=0$, we have
\[
\lim_{N \to \infty}(\nu\otimes\eta)
\big(
\left\{
(x,v)\in X\times\mathbb P\mathbb R^m \mid
S_N(x,v)\in(a,b)
\right\}
\big)
=
\Pb(S\in(a,b)).
\]
\end{theorem}

When the base system is assumed to be mixing, we obtain an affirmative answer to \textbf{Question~2}, yielding a two-sided mixing limit theorem:

\begin{theorem}[Mixing limit theorem for cocycles]
\label{thm:main}
Let $(X, \B, \mu, T)$ be an invertible mixing probability measure-preserving system, and let $C: X \times \Z \to \GL(m,\R)$ be an $m$-dimensional log-integrable cocycle over $T$. Suppose that $C$ satisfies a DLT on $(X\times \mathbb P\R^m,\mu\otimes \eta)$ with averaging sequence $\A = (A_N)$, normalizing sequence $\V = (V_N)$, and limiting distribution $S$. Let $\varphi_1, \varphi_2:X \to \R$ be two non-negative $L^{\infty}$ functions satisfying $\int_X \varphi_i\,\mathrm d\mu=1$ for $i=1,2$, and define a sequence of measures $\nu_N$ on $(X, \B)$ by 
\[
\nu_N(U)=\int_U \varphi_1\cdot(\varphi_2\circ T^N)\,\mathrm d\mu,\quad \forall U\in \B.
\]
Then the random variables $S_N$ on the probability spaces $\big(X\times \mathbb P\R^m,\frac{1}{\nu_N(X)}\nu_N \otimes\eta\big)$ converge in distribution to $S$ as $N \to \infty$.
\end{theorem}
   
In particular, given measurable sets $U_1, U_2 \in \B$ with $\mu(U_1), \mu(U_2) > 0$, we may specialize Theorem~\ref{thm:main} to the normalized indicator functions
\[
\varphi_i=\frac{1}{\mu(U_i)}\mathbf 1_{U_i},
\qquad i=1,2.
\]
Then, for every interval $(a,b) \subseteq \R$ satisfying $\Pb(S\in\{a,b\})=0$, we obtain the conditional limit law
\[
\lim_{N\to\infty} \frac{
(\mu\otimes\eta)
\left(
\left\{
(x,v) \mid
S_N(x,v)\in(a,b),\;
x\in U_1,\;
T^Nx\in U_2
\right\}
\right)
}{
\mu(U_1\cap T^{-N}U_2)
}
=
\Pb(S\in(a,b)).
\]
Since $T$ is mixing, we have $\mu(U_1\cap T^{-N}U_2) = \mu(U_1)\mu(U_2)+o(1)$ as $N \to \infty$, which immediately yields the asymptotic decoupling
\[
\lim_{N\to\infty} (\mu\otimes\eta)
\left(
\left\{
(x,v) \mid
S_N(x,v)\in(a,b),\;
x\in U_1,\;
T^Nx\in U_2
\right\}
\right)
=
\mu(U_1)\mu(U_2)\,
\Pb(S\in(a,b)).
\]

\section{Proof of the main results}

We begin with expressing the distributional convergence
through the following classical lemma.

\begin{lemma}\label{beginninglemma}
A sequence of real random variables $Z_n$ (not necessarily defined on
the same probability space) converges in distribution to a random
variable $Z$ if and only if, for every bounded Lipschitz function
$\alpha:\mathbb R\to\mathbb R$,
\[
    \E[\alpha(Z_n)]
    \longrightarrow
    \E[\alpha(Z)].
\]
\end{lemma}

The proof proceeds in four steps. First, we remove the dependence on
the projective variable by replacing the logarithmic expansion along a
direction with the logarithm of the operator norm. Second, we show
that the resulting observables on the base are asymptotically
invariant under every fixed time shift. This allows us to apply the
Eagleson averaging argument and remove the conditioning at the
initial time. Finally, a mixing estimate due to Gou\"ezel allows us
to handle the conditioning at the terminal time.

We begin with the reduction from the projective bundle to the base.
The key input is the following elementary linear algebra estimate.

\begin{lemma}\label{lem:bddaverage}
There exists a constant $\kappa_m \ge 0$, depending only on $m$, such that for every $A\in \mathrm{GL}(m,\mathbb R)$,
\[
\int_{\mathbb P\mathbb R^m}
\left|
\log \frac{\|Av\|}{\|v\|}
-
\log \opnorm{A}
\right|
\,\mathrm{d}\eta(v)
\le \kappa_m.
\]
\end{lemma}

\begin{proof}
For $m=1$, the projective space consists of a single point and the integrand vanishes identically, so we may take $\kappa_1 = 0$. Assume henceforth that $m \ge 2$. Throughout the proof, we identify each projective direction $v \in \mathbb{P}\mathbb{R}^m$ with a unit vector representative in $S^{m-1}$ (unique up to sign).

Let $e_A\in S^{m-1}$ be a unit right singular vector associated with the largest singular value of $A$, so that $\|Ae_A\|=\opnorm{A}$. Expanding any unit vector $v\in S^{m-1}$ in an orthonormal basis of right singular vectors of $A$ containing $e_A$ yields
\[
\|Av\|
\ge
\opnorm{A}\,|\langle v,e_A\rangle|.
\]
Since we also have $\|Av\|\le \opnorm{A}$, it follows that for every $v \in S^{m-1} \setminus e_A^\perp$ (and hence $\eta$-almost everywhere),
\[
0
\le
\log\opnorm{A} -\log\|Av\|
\le
-\log |\langle v,e_A\rangle|.
\]
Integrating over $\mathbb{P}\mathbb{R}^m$ gives
\[
\int_{\mathbb P\mathbb R^m}
\left|
\log \frac{\|Av\|}{\|v\|}
-
\log \opnorm{A}
\right|
\,\mathrm{d}\eta(v)
\le
\int_{\mathbb P\mathbb R^m}
-\log |\langle v,e_A\rangle|
\,\mathrm{d}\eta(v).
\]
By the $\mathrm{O}(m)$-invariance of $\eta$, the right-hand side is independent of the unit vector $e_A$. Thus, fixing an arbitrary unit vector $e\in S^{m-1}$ and setting
\[
\kappa_m
:=
\int_{\mathbb P\mathbb R^m}
-\log |\langle v,e\rangle|
\,\mathrm{d}\eta(v),
\]
it remains to verify that $\kappa_m<\infty$.

Let $\lambda_{m-1}$ denote the normalized $\mathrm{O}(m)$-invariant probability measure on $S^{m-1}$ and let $\pi:S^{m-1}\to\mathbb P\mathbb R^m$ be the canonical projection. Since $\eta=\pi_*\lambda_{m-1}$ and the integrand is even, pulling back to $S^{m-1}$ and using rotational invariance to choose $e=e_1$ gives
\[
\kappa_m
=
\int_{S^{m-1}}
-\log |\langle v,e\rangle|\,\mathrm{d}\lambda_{m-1}(v)
=
\int_{S^{m-1}}
-\log |\pi_1(v)|\,\mathrm{d}\lambda_{m-1}(v),
\]
where $\pi_1:S^{m-1}\to[-1,1]$, $\pi_1(v)=v_1$, is the projection onto the first coordinate. By a standard application of the coarea formula (see \cite[Chapter~3]{EvansGariepy}), the density $\rho_m$ of the push-forward measure $(\pi_1)_*\lambda_{m-1}$ on $(-1,1)$ is given by
\[
\rho_m(t)
=
c_m(1-t^2)^{\frac{m-3}{2}},
\qquad \text{where }
c_m
=
\frac{\Gamma(m/2)}
{\sqrt{\pi}\,\Gamma((m-1)/2)}.
\]
Consequently,
\[
\kappa_m=
\int_{-1}^1(-\log |t|)\,\rho_m(t)\,\mathrm{d}t=
2c_m
\int_0^1
(-\log t)(1-t^2)^{\frac{m-3}{2}}\,\mathrm{d}t.
\]
Since $-\log t$ is integrable near $t=0$ and $(1-t^2)^{(m-3)/2}$ is integrable near $t=1$ for all $m \ge 2$, we conclude that $\kappa_m < \infty$.
\end{proof}

Define
\[
S^{\mathrm{op}}_N(x)
\coloneqq
\frac{\sigma(x,N)-A_N}{V_N}=\frac{\log \opnorm{C(x,N)}-A_N}{V_N}.
\]
The previous lemma immediately gives the following:

\begin{proposition}
\label{prop:operator-reduction}
 $S_N$ converges in distribution to $S$ on
$(X\times\mathbb P\mathbb R^m,\mu\otimes\eta)$ if and only if
$ S_N^{\mathrm{op}}$ converges in distribution to $S$ on $(X,\mu)$.
\end{proposition}

\begin{proof}
By Lemma~\ref{beginninglemma}, it is enough to show that for every
bounded Lipschitz function $\alpha:\mathbb R\to\mathbb R$,
\[
\left|
\int_{X\times\mathbb P\mathbb R^m}
\alpha(S_N)\,\mathrm{d}(\mu\otimes\eta)
-
\int_X
\alpha(S_N^{\mathrm{op}})\,\mathrm{d}\mu
\right|
\longrightarrow 0
\qquad\text{as }N\to\infty.
\]

For $x\in X$, set
\[
F_N^\alpha(x)
:=
\int_{\mathbb P\mathbb R^m}
\alpha(S_N(x,v))\,\mathrm{d}\eta(v).
\]
Applying Lemma~\ref{lem:bddaverage} to $A=C(x,N)$ and using the
Lipschitz continuity of $\alpha$, we obtain
\[
\begin{aligned}
\left|
F_N^\alpha(x)
-
\alpha(S_N^{\mathrm{op}}(x))
\right|
&\le
\frac{\operatorname{Lip}(\alpha)}{V_N}
\int_{\mathbb P\mathbb R^m}
\left|
\sigma(x,v,N)-\sigma(x,N)
\right|\,\mathrm{d}\eta(v)
\\
&\le
\frac{\kappa_m\operatorname{Lip}(\alpha)}{V_N}.
\end{aligned}
\]
Since the bound is uniform in $x$, integrating over $X$ gives
\[
\left|
\int_{X\times\mathbb P\mathbb R^m}
\alpha(S_N)\,\mathrm{d}(\mu\otimes\eta)
-
\int_X
\alpha(S_N^{\mathrm{op}})\,\mathrm{d}\mu
\right|
\le
\frac{\kappa_m\operatorname{Lip}(\alpha)}{V_N}\longrightarrow 0.
\]

\end{proof}

Proposition~\ref{prop:operator-reduction} reduces the problem to a
distributional limit theorem for a sequence of real-valued functions
on the base space $X$. 

The next step is to establish the analogue of the finite-shift
invariance used in the proof of Eagleson's theorem. For Birkhoff sums, shifting the
base point by a fixed number of iterates changes only finitely many
terms at the two ends of the orbit segment. For logarithmic cocycle
growth, the same principle is encoded in the cocycle identity.

\begin{lemma}\label{lem:finite-shift}
Let $(X,\mathscr B,\mu,T)$ be an invertible probability
measure-preserving system, and let
\[
C:X\times\mathbb Z\longrightarrow\mathrm{GL}(m,\mathbb R)
\]
be a log-integrable cocycle. Fix a bounded Lipschitz function
$\alpha:\mathbb R\to\mathbb R$. For every fixed $k\in\mathbb N$,
for $1\leq j\leq k$,
\[
\left\|
\alpha(S_N^{\mathrm{op}})\circ T^j
-
\alpha(S_N^{\mathrm{op}})
\right\|_{L^1}
=
O_{\alpha,k}(V_N^{-1})
\]
as $N\to\infty$.
\end{lemma}

\begin{proof}
 Define
\[
q_j(x)
:=
\log^+\opnorm{C(x,j)}
+
\log^+\opnorm{C(x,j)^{-1}}.
\]
We first observe that $q_j\in L^1(X,\mu)$. Indeed, by the cocycle
identity,
\begin{align*}
    C(x,j) &=C(T^{j-1}x,1)\cdots C(Tx,1)C(x,1),\\
    C(x,j)^{-1}&=C(Tx,-1)C(T^2x,-1)\cdots C(T^{j}x,-1).
\end{align*}
and hence
\[
    q_j(x) 
\le
\sum_{\ell=0}^{j-1}
\log^+\opnorm{C(T^\ell x,1)}+\sum_{\ell=1}^{j}
\log^+\opnorm{C(T^\ell x,-1)}.
\]
Since $C$ is log-integrable and $T$ preserves $\mu$, the right-hand side
is integrable and hence $q_j\in L^1(X,\mu)$.

Notice that,
\[
C(T^jx,N)C(x,j)
=
C(x,N+j)
=
C(T^Nx,j)C(x,N),
\]
and hence
\[
C(T^jx,N)
=
C(T^Nx,j)\,C(x,N)\,C(x,j)^{-1}.
\]
By submultiplicativity of the operator norm, and by applying the same
argument to the inverse relation, we obtain
\begin{align*}
    \left|
\log\opnorm{C(T^jx,N)}
-
\log\opnorm{C(x,N)}
\right|&= \left| \log \frac{\opnorm{C(T^Nx,j)\,C(x,N)\,C(x,j)^{-1}}}{\opnorm{C(x,N)}}\right|\\
&\le
q_j(x)+q_j(T^Nx).
\end{align*}

Therefore, by the Lipschitz property of $\alpha$,
\[
\begin{aligned}
\left\|
\alpha(S_N^{\mathrm{op}})\circ T^j
-
\alpha(S_N^{\mathrm{op}})
\right\|_{L^1}
&\le
\frac{\operatorname{Lip}(\alpha)}{V_N}
\int_X
\bigl(q_j(x)+q_j(T^Nx)\bigr)\,\mathrm{d}\mu(x)
\\
&=
\frac{2\operatorname{Lip}(\alpha)\|q_j\|_{L^1}}{V_N},
\end{aligned}
\]
where we used the $T$-invariance of $\mu$ in the last equality.

Since $k$ is fixed,
\[
M_k:=\max_{1\leq j\leq k}\|q_j\|_{L^1}<\infty.
\]
Thus, uniformly for $1\leq j\leq k$,
\[
\left\|
\alpha(S_N^{\mathrm{op}})\circ T^j
-
\alpha(S_N^{\mathrm{op}})
\right\|_{L^1}
\le
\frac{2\operatorname{Lip}(\alpha)M_k}{V_N}
=
O_{\alpha,k}(V_N^{-1}),
\]
as claimed.
\end{proof}

At this point, all the cocycle-specific input has been reduced to the
finite-shift invariance of Lemma~\ref{lem:finite-shift}. The remaining
arguments are purely measure-theoretic and follow the same averaging
mechanism as in the scalar case.

Now we are ready to prove the theorem~\ref{thm:eagleson-cocycle}. 
\begin{proof}[Proof of theorem~\ref{thm:eagleson-cocycle}]
By Proposition~\ref{prop:operator-reduction}, it is enough to work with $ S_N^{\mathrm{op}}$ on the base. The proof then follows the usual
Eagleson averaging argument, with
Lemma~\ref{lem:finite-shift} providing the required asymptotic invariance under finite shifts. We include a proof here for completeness. 

Fix a bounded Lipschitz function $\alpha:\mathbb R\to\mathbb R$.
By Proposition~\ref{prop:operator-reduction}, it is enough to prove

\[
\int_X \alpha( S_N^{\mathrm{op}})\,\mathrm{d}\nu
\longrightarrow
\E[\alpha(S)].
\]

Let
\[
\rho=\frac{\mathrm{d}\nu}{\mathrm{d}\mu}\in L^1(X,\mu),
\qquad
\rho_0:=\rho-1\in L^1_0(X,\mu).
\]
Then
\[
\int_X \alpha(   S_N^{\mathrm{op}})\,\mathrm{d}\nu
=
\int_X \alpha(   S_N^{\mathrm{op}})\,\mathrm{d}\mu
+
\int_X \alpha(   S_N^{\mathrm{op}})\rho_0\,\mathrm{d}\mu.
\]
By the assumed distributional convergence under $\mu$, the first
term converges to $\mathbb E[\alpha(S)]$. It remains to prove that
\[
\int_X \alpha(   S_N^{\mathrm{op}})\rho_0\,\mathrm{d}\mu
\longrightarrow0.
\]

Fix $\varepsilon>0$. Choose
$\widetilde\rho_0\in L^\infty_0(X,\mu)$ such that
\[
\|\rho_0-\widetilde\rho_0\|_{L^1}<\varepsilon.
\]
Since $T$ is ergodic and $\widetilde\rho_0$ has zero mean, the mean
ergodic theorem gives a $k\in\mathbb N$ such that
\[
\left\|
\frac1k\sum_{j=0}^{k-1}\widetilde\rho_0\circ T^j
\right\|_{L^2}
<\varepsilon.
\]

For this fixed $k$, the $T$-invariance of \(\mu\) and Lemma~\ref{lem:finite-shift} implies, for each
$j=0,\ldots,k-1$,
\begin{align*}
    \int_X
\alpha(   S_N^{\mathrm{op}})\widetilde\rho_0\,\mathrm{d}\mu
&= \int_X \alpha(S_N^{\mathrm{op}})\circ T^j \cdot(\widetilde \rho_0 \circ T^j)\,\mathrm{d}\mu \\
&=\int_X
\big(\alpha(S_N^{\mathrm{op}})\circ T^j-\alpha(   S_N^{\mathrm{op}})\big)
(\widetilde\rho_0\circ T^j)\,\mathrm{d}\mu
+ \int_X \alpha(S_N^{\mathrm{op}})(\widetilde\rho_0\circ T^j)\,\mathrm{d}\mu\\
&= \int_X \alpha(S_N^{\mathrm{op}})(\widetilde\rho_0\circ T^j)+o(1),
\end{align*}
as $N\to\infty$. Averaging over $j$ yields
\[
\int_X
\alpha(   S_N^{\mathrm{op}})\widetilde\rho_0\,\mathrm{d}\mu
=
\int_X
\alpha(   S_N^{\mathrm{op}})
\left(
\frac1k\sum_{j=0}^{k-1}\widetilde\rho_0\circ T^j
\right)
d\mu
+
o(1).
\]
Hence, by Cauchy--Schwarz,
\[
\begin{aligned}
\left|
\int_X \alpha(   S_N^{\mathrm{op}})\rho_0\,\mathrm{d}\mu
\right|
&\le
\|\alpha\|_{L^\infty}
\|\rho_0-\widetilde\rho_0\|_{L^1}
\\
&\quad+
\left\|\alpha(   S_N^{\mathrm{op}})\right\|_{L^2}\cdot
\left\|
\frac1k\sum_{j=0}^{k-1}\widetilde\rho_0\circ T^j
\right\|_{L^2}
+o(1)
\\
&\le
2\varepsilon\|\alpha\|_{L^\infty}+o(1).
\end{aligned}
\]
It follows that
\[
\limsup_{N\to\infty}
\left|
\int_X \alpha(   S_N^{\mathrm{op}})\rho_0\,\mathrm{d}\mu
\right|
\le
2\varepsilon\|\alpha\|_{L^\infty}.
\]
Since $\varepsilon>0$ is arbitrary,
\[
\int_X \alpha(   S_N^{\mathrm{op}})\rho_0\,\mathrm{d}\mu\to0.
\]
Therefore
\[
\int_X \alpha(   S_N^{\mathrm{op}})\,\mathrm{d}\nu
\longrightarrow
\mathbb \E[\alpha(S)].
\]
By Lemma~\ref{beginninglemma}, this proves the desired convergence in
distribution under $\nu$, and Proposition~\ref{prop:operator-reduction}
transfers the conclusion back to the projective cocycle variables
$S_N(x,v)$.
\end{proof}

Theorem~\ref{thm:eagleson-cocycle} shows that conditioning on the
initial position does not alter the limiting distribution. The only
remaining difficulty is the conditioning at time $N$, whose location
moves with $N$. To deal with this, we need the following lemma due to Gou\"ezel about the mixing systems:

\begin{lemma}
\label{lem:gouezel-averaging}
Let $(X,\mathscr B,\mu,T)$ be a mixing probability measure-preserving
system, $\psi_1,\psi_2\in L^{\infty}(X,\mu)$ with $\int_X \psi_2\ d\mu = 0$. Then 
    $$ \limsup_{N \to \infty}\left\| \frac{1}{k} \sum_{j=0}^{k-1} \psi_1 \circ T^j \cdot \psi_2 \circ T^{N+j} \right\|_{L^2} =o(1), \quad \text{as } k \to \infty$$
\end{lemma}

\begin{proof}
See Gou\"ezel~\cite[Lemma~2.3]{GOUËZEL_2020}.
\end{proof}

\begin{proof}[Proof of Theorem~\ref{thm:main}]
    We are now ready to prove the main theorem. By
Proposition~\ref{prop:operator-reduction}, it is enough to work with
$S_N^{\mathrm{op}}$. Let $\alpha:\mathbb R\to\mathbb R$ be bounded and
Lipschitz. Our goal is to show that 
    \begin{equation}\label{MDLTgoal}
        \int_X \alpha(S_N^{\mathrm{op}})\cdot(\varphi_1 \cdot \varphi_2\circ T^N)\ d\mu \longrightarrow\left(\int_X\varphi_1\ d\mu\right) \E(\alpha(S)) \left(\int_X\varphi_2\ d\mu\right)
    \end{equation}
    Notice that when \(\varphi_2\) is constant, the convergence~(\ref{MDLTgoal}) holds by Theorem~\ref{thm:eagleson-cocycle}. Hence, we can without loss of generality replace \(\varphi_2\) with \(\varphi_2-\int \varphi_2d\mu\) and assume \(\int\varphi_2 d\mu=0\). It remains to show that
\[
I_N\coloneqq\int_X \alpha(S^{\mathrm{op}}_N)\cdot(\varphi_1\cdot\varphi_2\circ T^N)\,\mathrm{d}\mu
\longrightarrow 0.
\]

Fix $k\geq1$ and, for $0\leq j<k$, write $\Phi_{N,j}=\varphi_1\circ T^j \cdot \varphi_2\circ T^{N+j}$, then the $T$-invariance of $\mu$ gives
\[
I_N=\int_X\alpha(S_N^{\mathrm{op}}\circ T^j)\cdot \Phi_{N,j}\,\mathrm{d}\mu.
\]
Consequently,
\[
\left|I_N-\int_X \alpha(S_N^{\mathrm{op}}) \cdot\Phi_{N,j}\,\mathrm{d}\mu\right|
\leq \|\varphi_1\|_{L^\infty}\|\varphi_2\|_{L^\infty}
\|\alpha(S_N^{\mathrm{op}}\circ T^j)-\alpha(S_N^{\mathrm{op}})\|_{L^1}.
\]
For each fixed $j$, the right-hand side tends to zero by
Lemma~\ref{lem:finite-shift}; Averaging these
identities over $j=0,\ldots,k-1$, with $k$ fixed, yields
\[
I_N=\int_X \alpha(S_N^{\mathrm{op}})\cdot
\left(\frac1k\sum_{j=0}^{k-1}
\varphi_1\circ T^j\,\varphi_2\circ T^{N+j}\right)d\mu
+o_k(1).
\]
By Cauchy--Schwarz and the boundedness of $\alpha$,
\[
\limsup_{N\to\infty}|I_N|
\leq \|\alpha\|_{L^\infty}
\limsup_{N\to\infty}
\left\|\frac1k\sum_{j=0}^{k-1}
\varphi_1\circ T^j\,\varphi_2\circ T^{N+j}
\right\|_{L^2}.
\]
Lemma~\ref{lem:gouezel-averaging} makes the right-hand side tend to
zero as $k\to\infty$. Thus $I_N\to0$, and the proof is finished.
\end{proof}

\bibliographystyle{amsplain}  % 或 amsalpha, plain, abbrv, alpha 等
\bibliography{references}
\end{document}